\documentclass[12pt]{amsart}
\usepackage{amsmath,amssymb,amsthm}
\usepackage[]{latexsym}
\newtheorem{defn0}{Definition}[section]
\newtheorem{prop0}[defn0]{Proposition}
\newtheorem{thm0}[defn0]{Theorem}
\newtheorem{lemma0}[defn0]{Lemma}
\newtheorem{corollary0}[defn0]{Corollary}
\newtheorem{example0}[defn0]{Example}
\newtheorem{remark0}[defn0]{Remark}
\newtheorem{conjecture0}[defn0]{Conjecture}
\newtheorem{notation0}[defn0]{Notation}

\newenvironment{definition}{\begin{defn0}\rm}{\end{defn0}}
\newenvironment{proposition}{\begin{prop0}}{\end{prop0}}
\newenvironment{theorem}{\begin{thm0}}{\end{thm0}}

\newenvironment{corollary}{\begin{corollary0}}{\end{corollary0}}

\newenvironment{remark}{\begin{remark0}\rm}{\end{remark0}}

\newcommand{\Hom}{{\mathrm {Hom}}}

\newcommand{\Gal}{{\mathrm Gal}}
\newcommand{\n}{\mathrm{n}}
\newcommand{\tr}{\mathrm{tr}}

\newcommand{\Pic}{\mathrm{Pic}}
\newcommand{\Norm}{\mathrm{N}}
\newcommand{\Aut}{\mathrm{Aut}}

\newcommand{\End}{{\mathrm{End}}}

\newcommand{\disc}{{\mathrm disc}}
\newcommand{\NS}{{\mathrm NS}}
\newcommand{\M}{{\mathrm M}}

\newcommand{\Z}{{\mathbb Z}}
\newcommand{\Q}{{\mathbb Q}}
\newcommand{\C}{{\mathbb C}}
\newcommand{\R}{{\mathbb R}}

\newcommand{\cL}{{\mathcal L}}

\newcommand{\cI}{{\mathcal I}}
\newcommand{\cA}{{\mathcal A}}
\newcommand{\cQ}{{\mathcal Q}}
\newcommand{\cH}{{\mathcal H}}

\newcommand{\cN}{{\mathcal N}}

\newcommand{\cO}{{\mathcal O}}

\newcommand{\sk}{\vspace{0.1in}}

\newcommand{\ra}{\rightarrow}

\newcommand{\om}{{\omega }}

\include{thebibliography}

\begin{document}

\title{Modular Shimura varieties and forgetful maps}

\author{Victor Rotger}
\footnote{Partially supported by a grant FPI from the
Ministerio de Ciencia y Tecnolog\'{\i}a BFM2000-0627.}

\address{Universitat Polit\`{e}cnica de Catalunya,
Departament de Matem\`{a}tica Aplicada IV (EUPVG), Av.\
Victor Balaguer s/n, 08800 Vilanova i la Geltr\'{u}, Spain.}

\email{vrotger@mat.upc.es}

\subjclass{11G18, 14G35}

\keywords{Shimura variety, moduli space, abelian variety, quaternion
algebra}

$$
\mbox{Appeared in {\em Trans.\ Amer. Math.\ Soc.\ } {\bf 356} (2004), 1535-1550.}
$$

\begin{abstract}
In this note we consider several maps that occur naturally between
modular Shimura varieties, Hilbert-Blumenthal varieties and the
moduli spaces of polarized abelian varieties when forgetting
certain endomorphism structures. We prove that, up to birational
equivalences, these forgetful maps coincide with the natural
projection by suitable abelian groups of Atkin-Lehner involutions.
\end{abstract}

\maketitle

\section*{Introduction}
\sk\noindent

Let $F$ be a totally real number field of degree $[F:\Q ]=n$ and
let $B$ be a totally indefinite quaternion algebra over $F$. In
this note we will be concerned with certain Shimura varieties
$X_B$ attached to an arithmetic datum arising from $B$ and several
maps that occur naturally between them.

As complex manifolds, these varieties can be described as
quotients of certain bounded symmetric domains by arithmetic
groups acting on them and, by the theory of Baily-Borel
(\cite{BaBo}), they become quasi-projective complex algebraic
varieties. Shimura discovered a moduli interpretation of these
varieties which allowed him to construct canonical models $X_B/\Q
$ over the field $\Q $ of rational numbers. Shimura also explored
their arithmetic, showing that the coordinates of so-called {\em
Heegner points} on $X_B$ generate certain class fields and that
the Galois action on them can be described by explicit reciprocity
laws (\cite{Sh1}, \cite{Sh2}).

The nature of the Shimura varieties $X_B$ differs notably
depending on the existence or absence of zero divisors in $B$.
When $B$ is the split algebra $\mathrm {M}_2(F)$, the varieties
$X_B$ are classically called {\em Hilbert-Blumenthal modular
varieties}. These are non-complete and suitable compactifications
of them can be constructed, though at the cost of producing new
singularities. The literature on them is abundant, especially on
the low dimensional cases (in dimension $1$, these are called
modular curves and have become crucial in many aspects of number
theory; in dimension $2$, a reference to Hilbert modular surfaces
is van der Geer's book \cite{vdGe}).

On the other hand, when $B$ is non-split, that is, it is a
division algebra, then the emerging Shimura varieties $X_B/\Q $
are already projective. This fact makes the study of their
arithmetic highly difficult since, in the Hilbert-Blumenthal case,
much of it is encoded in the added cusps. In remarkable contrast
to Hilbert-Blumenthal varieties, Shimura proved that, when $B\not
\simeq \M _2(F)$, the varieties $X_B$ do not have real points and
therefore do not have rational points over any number field that
admits a real embedding. In the last years, there has been
increasing interest on Shimura curves arising from rational
indefinite quaternion algebras, since they play a crucial role in
modularity questions (cf.\,\cite{HaHaMo} or \cite{Ri}, for
instance).

From the modular point of view, there are natural maps
$$
\pi : X_B\rightarrow \mathcal A _{g, (d_1, ..., d_g)}
$$
from $X_B$ into the moduli spaces of polarized abelian varieties
that involve {\em forgetting} some additional structures. Further,
for any totally real field $L$ containing $F$ and embedded in $B$,
these maps factor through a morphism $\pi _L: X_B\rightarrow
\mathcal H _L$ into a Hilbert-Blumenthal modular variety
(cf.\,Section 2).

It is the purpose of this note to describe in detail the nature of
these morphisms and their image in the several moduli spaces of
abelian varieties. As we will see, $\pi $ and $\pi _L$ are
quasifinite maps and the mere computation of their degree turns
out to be unexpectedly subtle, as it is based on arithmetic
questions on $B$ that were recently studied by Chinburg and
Friedman in a number of papers in relation with arithmetic
hyperbolic 3-orbifolds (cf. \cite{ChFr1}, \cite{ChFr2}).

In the first section we briefly recall some basic facts on
quaternion algebras that we need throughout the article. Next, in
Section 2 we review the construction due to Shimura of the
varieties we will be dealing with and the above mentioned maps. In
the third section, we introduce several Atkin-Lehner groups of
automorphisms acting on them and describe their modular
interpretation. In Section 4 we introduce what we call the {\em
stable} and {\em twisting Atkin-Lehner groups} and we state our
main result in Theorem \ref{main}).

The rest of the note is devoted to presenting the proof of Theorem
\ref{main}. To do so, in Section \ref{Bpart} we introduce what we
call the $B$-part of the N\'{e}ron-Severi group of an abelian variety
whose algebra of endomorphisms $\End (A)\otimes \Q $ contains the
quaternion algebra $B$.

In closing this introduction, we would point out two different
applications of the results presented in this note. The first one
concerns the geometry of the quaternionic locus $\cQ _{\cO }=\{
[A, \cL ], \End (A)\supseteq \cO \} $ of abelian varieties
admitting multiplication by a maximal order $\cO $ in the moduli
space $\cA _g$ of principally polarized abelian varieties of even
dimension $g$. By means of Theorem \ref{main} and Eichler's theory
on optimal embeddings, the number of irreducible components of
$\cQ _{\cO }$ can be related to certain class numbers and
conditions can be given for its irreducibility. The emerging
picture bears a remarkable resemblance to the classical theory of
genus of quadratic forms and {\em idonei numbers}
(cf.\,\cite{Ro4}).

Secondly, Theorem \ref{main} can also be used to explore the
arithmetic of QM-abelian varieties. Indeed, in \cite{Ro3}, we use
it to describe the field of moduli of the quaternionic
endomorphisms acting on them.

Finally, the results of this paper are also used by Clark in
\cite{Cl} in relation to abelian surfaces with potential
quaternionic multiplication and their moduli.

{\em Acknowledgements. } I am indebted to P. Bayer for her
assistance throughout the elaboration of this work. I also express my
gratitude to B. Mazur and P. Clark for the interest they have shown in
this article. Finally, I thank J. Kramer
and U. Kuehn for their warm hospitality at the
Humboldt-Universit\"{a}t zu Berlin during the fall of 2001.

\section{An overview of quaternion algebras}

In this section let $F$ denote an arbitrary number field. Let $B$
be a quaternion algebra over $F$, that is, a simple algebra of
rank$_F (B) = 4$ over its centre $F$. Elements $i$, $j\in B$ can
be chosen such that $B = F + F i + F j + F i j$ with the relations
$i j = -j i $ and $i^2=a$, $j^2=b\in F^*$. The classical notation
for it is $B = \left ( \frac {a, b}{F}\right )$.

A place $v$ of $F$, archimedean or not, is said to ramify in $B$ if
$B\otimes F_v$ is a division algebra over the completion $F_v$ of $F$
at $v$. When $v$ is archimedean, in this case we also say that $B$ is
{\em definite} at $v$ and otherwise we call it {\em indefinite}.
A quaternion algebra that is indefinite at all archimedean places
is called {\em totally indefinite}.

There are finitely many places of $F$ that ramify in a quaternion
algebra $B$ and they occur in an even number. The discriminant ideal of
$B$ is the (square-free) product $\mathrm{disc }(B) = \prod \mathfrak p$ of
the finite primes that ramify in $B$.

Let $\beta = x+y i+z j+t i j\mapsto \bar \beta =x-y i-z j-t i j$
denote the conjugation map on $B$. Then the reduced trace and
reduced norm of $\beta $ are defined as $\mbox{tr }(\beta )=\beta
+\bar \beta \in F$ and $\mbox{n }(\beta )=\beta \bar \beta \in F$
respectively. An element $\beta \in B$ of null reduced trace $\tr
(\beta ) = 0$ is called {\em pure} and, for any sub-module $\cI $
of $B$, we will denote by $\cI _0\subseteq B_0$ the group of pure
quaternions of $\cI $. An element $\beta \in B$ is integral over
the ring of integers $R_F$ of $F$ if both $\tr (\beta )$ and $\n
(\beta )$ lie in $R_F$. An order $\cO $ of integers in $B$ is an
$R_F$-finitely generated subring of integral elements such that
$\cO \otimes _{R_F}F=B$. It is called a {\em maximal order} if it
is not properly contained in any other. Unlike number fields,
maximal orders are not unique. The {\em type number} $t(B)$ of $B$
is the number of $B^*$-conjugation classes of maximal orders in
$B$.

A left ideal of a maximal order $\cO $ is an $R_F$-finitely
generated module $\cI \subset B$ with $\cI \otimes _{R_F} F = B$
and $\cO \cdot \cI = \cI $. The ideal is two-sided if in addition
$\cI \cdot \cO = \cI $. We will often consider $\cI $ in
$\mbox{Pic}_{\ell }(\cO )$, that is, up to principal $\cO $-left
ideals $\cO \beta $, $\beta \in B^*$. See \cite{Vi}, p.\,25, for
details. The conjugate ideal of an ideal $\cI $ in $B$ is the
ideal $\bar {\cI } = \{ \bar {\beta }, \beta \in \cI \} $. If $\cI
$ is a left ideal of a maximal order $\cO $, then $\bar {\cI }$ is
a right $\cO $-ideal. The product ideal $\cI \cdot \bar {\cI }$ is
then a two-sided ideal of $\cO $ and we actually have that $\cN
(\cI ):=\n (\cI )\cO =\cI \cdot \bar {\cI }$.

Finally, if $K\subset F$ is any sub-field of $F$, the {\em
codifferent } of $\cI $ over $K$ is defined to be $\cI ^{\sharp
}_{B/K}=\{\beta \in B: \mathrm{tr}_{B/K}(\cI \beta )\subseteq
R_K\}$. It is a right ideal of $\cO $. In our considerations on
N\'{e}ron-Severi groups, the group $\cN (\cI )^{\sharp }_0$ of pure
quaternions of the codifferent of the norm ideal  $\cN (\cI )$ over
$\Q $ will play an
important role (cf.\,Section \ref{Bpart}).

Assume now that $F$ is totally real and let $F^*_+$ denote the
group of totally positive elements. We will also let $R^*_{F_+}$
stand for the group of totally positive units in $R_F^*$. A
positive (anti-)involution $\varrho $ on $B$ is a map $\varrho
:B\rightarrow B$ such that $(\beta _1+\beta _2)^{\varrho }=\beta
_1^{\varrho }+\beta _2^{\varrho } $ and $ (\beta _1\cdot \beta
_2)^{\varrho }=\beta _2^{\varrho }\cdot \beta _1^{\varrho } $ for
any $\beta _1$, $\beta _2\in B$, and such that $\tr (\beta \cdot
\beta ^{\varrho })\in F^*_+ $ for any $\beta \in B^*$.

By the Skolem-Noether Theorem, if $\varrho $ is a positive
(anti-)involution, there exists $\mu \in B^*$ such that $\beta
^{\varrho }=\mu ^{-1}\bar \beta \mu$. Further, it is easily shown
(\cite{Ro2}) that the positiveness of $\varrho $ implies that $\tr
(\mu )=0$ and $\n (\mu )\in F^*_+$. The element $\mu $ is
determined up to multiplication by elements of $F^*$ and we will
sometimes use the notation $\varrho = \varrho _{\mu }$.

\section{Shimura varieties and maps between them}
\label{ShimuraVarieties}

Let $F$ be a totally real number field of degree $n$ and $B$ a totally
indefinite quaternion algebra over it. Call $\mathcal D = \mathrm{disc}(B) =
\mathfrak p_1\cdot ...\cdot \mathfrak p_{2r}$ its discriminant.

Fix a datum $( \cO , \cI , \varrho )$ consisting of a maximal
order $\cO $ in $B$, a left $\cO $-ideal $\cI $ (or rather its
class in $\mathrm{Pic}_{\ell }(\cO )$) and a positive involution
$\varrho _{\mu }$ with $\mu \in \cO $, $\mu ^2+\delta =0$, $\delta
\in F^*_+$. Attached to the datum $( \cO , \cI , \varrho )$ there
is the following moduli problem over $\Q $: classifying
isomorphism classes of triplets $(A, \iota , \cL )$ where

\begin{itemize}
\item

$A$ is an abelian variety of dimension $g=2n$.

\item

$\iota :\cO \hookrightarrow \mbox{ End }(A) $ is a ring homomorphism
such that $H_1(A, \Z )$, regarded as a left $\cO $-module, is
isomorphic to the left ideal $\cI $.

\item
$\cL $ is a pri\-mitive pola\-riza\-tion on $A$ such that the Rosa\-ti involu\-tion
$\circ :\mathrm{End }^0(A)\\ \rightarrow \mathrm{End }^0(A)$ with
respect to $\cL $ on $\mathrm{End }^0(A) = \mathrm{End }(A)\otimes _{\Z }\Q $
coincides with $\varrho $ when restricted to $\iota (\cO )$: $\circ
|_{\iota (\cO )}=\varrho \cdot \iota $.

\end{itemize}

This triplet will be referred to as a polarized abelian variety
with multiplication by $\cO $. Two triplets $(A_1, \iota _1, \cL
_1)$, $(A_2, \iota _2, \cL _2)$ are isomorphic if there exists an
isomorphism $\alpha \in \Hom (A_1, A_2)$ such that $\alpha \iota
_1 (\beta ) = \iota _2(\beta ) \alpha $ for any $\beta \in \cO $
and $\alpha ^* (\cL _2) = \cL _1\in \NS (A_1)$. We recall that a
polarization $\cL \in \NS (A)$ is called primitive if $\cL \not
\in d\cdot \NS (A)$ for any $d\in \Z, d\geq 2$. Its type is then
$(1, d_2, ..., d_g)$ for $d_i\mid d_{i+1}$, $i=2, ..., g-1$.
Remark also that, since a priori there is not a canonical
structure of $R_F$-algebra on $\End (A)$, the immersion $\iota
:\cO \hookrightarrow \End (A)$ is just a homomorphism of rings (or
$\Z $-algebras).

The corresponding moduli functor is coarsely represented by an
irreducible and reduced quasi-projective scheme $X_B/\Q = X_{(\cO
, \cI , \varrho )}/\Q $ over $\Q $ and of dimension $n=[F: \Q ]$.
Moreover, if $B$ is division (that is, $r>0$), the Shimura variety
$X_B$ is complete (cf.\,\cite{Sh1}, \cite{Sh2}).

Complex analytically, the manifold $X_B(\C )$ can be described as
the quotient of a symmetric space by the action of a discontinuous
group as follows. Since $B$ is totally indefinite, we may fix an
embedding $B\hookrightarrow B\otimes _{\Q }\R \stackrel{\simeq
}{\rightarrow }\mbox{M}_2(\R )\oplus ...\oplus \M _2(\R )$ and
regard the group $\Gamma _B = \cO ^1 = \{\gamma \in \cO ^*,
\mbox{n}(\gamma )=1\} $ as a discrete subgroup of $\mbox{SL}_2(\R
)^n$. An element $\gamma = (\gamma _1, ..., \gamma _n)\in \Gamma
_B$ acts on the cartesian product ${\bf \mathfrak H}^n$ of $n$
copies of Poincar\'{e}'s upper half plane ${\bf \mathfrak H}= \{ x+y
i, x, y\in \R , y>0\} $ by Moebius transformations: $\gamma \cdot
(\tau _1, ..., \tau _n)^t = (\frac{a_1 \tau _1+b_1}{c_1 \tau
_1+d_1}, ..., \frac{a_n \tau _n+b_n}{c_n \tau _n+d_n})^t$ where
$\gamma _i = \begin{pmatrix}
  a_i & b_i \\
  c_i & d_i
\end{pmatrix}\in \mathrm{SL}_2(\R )$. Then

$$
\Gamma _B\backslash \mathfrak H^n\simeq X_B(\C ).
$$

\subsection{The maps into $\cA _g$, $\cH _F$ and $\cH _L$}
\label{Shimuramaps}

Let now $(A, \iota ,\cL )$ be a (primitively) polarized abelian
variety with multiplication by $\cO $ with respect to the datum
$(\cO , \cI ,\varrho )$. Note that the type $(1, d_2, ..., d_g)$
of $\cL $ was not specified when posing the above moduli problem.
However, since $X_{(\cO , \cI , \varrho )}$ is connected, the type
of $\cL $ only depends on the datum $(\cO , \cI , \varrho )$ and
not on the particular triplet $(A, \iota , \cL )$. Thus $(1, d_2,
..., d_g)$ will often be referred to as the type of $(\cO , \cI ,
\varrho )$. This observation allows us to define a natural
morphism

$$
\begin{matrix}
 \pi : & X_{(\cO , \cI , \varrho )} &\longrightarrow  & \mathcal A _{g,
(1, d_2, ..., d_g)}  \\
       &(A, \iota , \cL )& \mapsto & (A, \cL )
\end{matrix}
$$
from the Shimura variety to the moduli space of polarized abelian
varieties of type $(1, d_2, ..., d_g)$ that consists of {\em
forgetting} the quaternionic endomorphism structure. This morphism
is representable, proper and defined over $\Q$. Moreover, as we
now explain, the morphism $\pi :X\rightarrow \mathcal A _{g, (1,
d_2, ..., d_g)}$ factors in a natural way through certain Hilbert
modular varieties.

\begin{definition}

An Eichler pair $(S, \varphi )$ for $\cO $ is a pair consisting of

\begin{itemize}

\item
An order $S$ over $R_F$ in a quadratic extension $L$ of $F$ and

\item An $R_F$-embedding $\varphi :S\hookrightarrow \cO $ such that
$\varphi (S)=\varphi (L)\cap \cO $.

\end{itemize}
An Eichler pair is totally real if $L$ is.

\end{definition}

Note that not all orders $S$ in quadratic extensions $L$ of $F$
can be embedded in $\cO $. Namely, by Eichler's theory on optimal
embeddings, (\cite{Vi}, 5. C ), there exists an embedding $\varphi
$ of the ring of integers $R_L$ of $L$ into $\cO $ iff any prime
ideal $\mathfrak p$ of $F$ that ramifies in $B$ either remains
inert or ramifies in $L$. Here, the fact that $B$ is division and
splits at least at one archimedean place makes the condition for
the embeddability of $R_L$ in $\cO $ particularly neat. Otherwise,
it depends heavily on the conjugation class of $\cO $.

Given a totally real Eichler pair $(S , \varphi )$, we
can then consider the Hil\-bert mo\-dular va\-rie\-ty $\mathcal H _S$ that
classifies isomorphism classes of triplets $(A, i, \cL )$ where

\begin{itemize}

\item
$A$ is an
abelian variety of dimension $[L:\Q ] = 2n$,

\item
$i:S\hookrightarrow \End
(A)$ is a ring homomorphism and

\item
$\cL $ is a polarization of type $(1, d_2, ..., d_g)$ on $A$.

\end{itemize}

The scheme $\mathcal H _S$ is $2 n$-dimensional, non-complete and
defined over $\Q $. Actually, this is a particular case of the
Shimura varieties mentioned above: $\cH _S$ is the union of
several irreducible components, all them isomorphic to Shimura
varieties $X_{(\mathrm{M}_2(S), \cI , \varrho )}$ for several
$\mathrm{M}_2(S)$-left ideals $\cI $ and involutions $\varrho $.
Imitating the construction of $\pi $ we then obtain, for any
totally real Eichler pair $(S, \varphi )$, a morphism

$$
\begin{matrix}
\pi _{(S, \varphi )}: & X_{(\cO , \cI , \varrho)} &
\longrightarrow & \mathcal H _S \\
& (A, \iota , \cL ) & \mapsto & (A, \iota \cdot \varphi , \cL )
\end{matrix}
$$
where $\iota \cdot \varphi :S\hookrightarrow \cO
\hookrightarrow \End (A)$. Finally, and in a
similar way, we also have a morphism

$$
\pi _F:\quad X_{(\cO , \cI , \varrho )}\longrightarrow \mathcal H
_F
$$
from $X_B$ into a Hilbert variety $\cH _F$. However, this time
$\cH _F$ does not correspond to any of the Shimura varieties $X_B$
introduced above for any quaternion algebra $B$ (not even
$\mathrm{M}_2(F)$). On this occasion, by $\cH _F$ we mean the
(reduced) scheme over $\Q $ that coarsely represents the functor
attached to the moduli problem of classifying $(1, d_2, ...,
d_g)$-polarized abelian varieties $A$ of dimension $2 n$ together
with an homomorphism $R_F\hookrightarrow \End (A)$. The variety
$\cH _F$ has dimension $3 n$ and $\cH _F(\C )$ is the quotient of
$n$ copies of the $2$-Siegel half space $\mathfrak H_2$ by a
suitable discontinuous group (cf.\,\cite{Sh1}, \cite{LaBi},
Chapter 9). Notice that, when $F=\Q $, $\cH _F = \cA _{2, (1, d)}$
is Igusa's three-fold of level $d\geq 1$, the moduli space of $(1,
d)$-polarized abelian surfaces.

Then, the map $\pi _F:\quad X_{(\cO , \cI , \varrho
)}\longrightarrow \cH _F$ is constructed as above: by restricting
the endomorphism structure of a triplet $(A, \iota , \cL )$ to
$\iota |_{R_F}:R_F\hookrightarrow \End (A)$. The whole picture can
be summarized with the following commutative diagram of morphisms
between Shimura varieties

$$
\begin{matrix}
   &    &    & \cH _{S_1} &    &    &   &   \\
 & &\stackrel {\pi _{(S_1, \varphi _1)}}{\nearrow }&\stackrel{:}{.}
 &\searrow & & & \\
\pi : & X_B & \stackrel {\pi _{(S_2, \varphi _2)}}
{\longrightarrow } & \cH _{S_2} & \longrightarrow & \cH _F &
\rightarrow & \mathcal A _{g, (1, d_2, ..., d_g)} \\
  &  & \stackrel {\searrow }{\pi _{(S_n, \varphi _n)}} & \stackrel{:}{.} &
\nearrow  & & & \\
  &  &  & \cH _{S_n} &  &  &  &
\end{matrix}
$$
Note that, while there is a canonical forgetful map from $X_B$ to
$\cH _F$, we obtain distinct maps from $X_B$ to $\cH _S$ as
$\varphi :S\hookrightarrow \cO $ varies among all possible Eichler
embeddings.

\section{The Atkin-Lehner group of a Shimura variety}
\label{Atkin}

As before, let $F$ be a totally real number field of degree $n$. We
denote by $\Pic (F)$ the class group of fractional ideals up to
principal ideals of $F$. Similarly, $\Pic _+(F)$ will stand for
the narrow class group of fractional ideals up to totally positive
principal ideals $a R_F$, $a\in F^*_+$, of $F$.

Let $B$ be a totally indefinite quaternion algebra over $F$ of
discriminant $\mathcal D = \mathfrak p_1\cdot ...\cdot \mathfrak
p_{2r}$ and let $\cO $ be a maximal order in $B$. Define the
groups $\cO ^*\supset \cO ^*_+\supseteq \cO ^1$ of units in $\cO
$, units in $\cO $ of totally positive reduced norm and units in
$\cO $ of reduced norm $1$ respectively. By the
Hasse-Schilling-Maass Theorem (cf.\,\cite{HaSc}, \cite{Vi},
p.\,90) we have the exact sequence

$$
1\rightarrow \cO ^1\rightarrow \cO ^*_+\stackrel{\n }{\rightarrow }
R^*_{F_+}\rightarrow 1.
$$

\begin{definition}

Let $\mathrm{N} _{B^*}(\cO ) = \{ \gamma \in B^*, \gamma ^{-1}\cO
\gamma \subseteq \cO \} $ be the normalizer of $\cO $ in $B^*$.
The {\em Atkin-Lehner group} $W$ of $\cO $ is

$$
W=\mathrm{N} _{B^*}(\cO )/ F^* \cO ^*.
$$
Let $B^*_+$ denote the group of invertible quaternions of $B$ of
totally positive reduced norm. We define the {\em positive Atkin-Lehner groups}

$$
W_+ = \mathrm{N} _{B_+^*}(\cO )/F^* \cO ^*_+
$$
and
$$
W^1 = \mathrm{N} _{B_+^*}(\cO )/F^* \cO ^1.
$$
\end{definition}

Let us note that, by the Skolem-Noether Theorem, $\mathrm{N}
_{B^*}(\cO )/F^*$ coincides with the group of automorphisms of
$\cO $. The Atkin-Lehner group $W$ is identified with the group of
principal two-sided ideals of $\cO $ by the assignation $\omega
\in W \mapsto \cO \cdot \omega $ and it is a finite abelian
$2$-group. More precisely, by Eichler's results (\cite{Ei1},
\cite{Ei2}, \cite{Vi}, Theorem 5.7), the reduced norm $\n
:B^*\rightarrow F^*$ induces an isomorphism $W\simeq \Z /2 \Z
\times \stackrel{t}{...} \times \Z /2 \Z $ for some $t\le 2 r$.
Moreover, any involution $[\om ]$ on $W$ can be represented by an
element $\om \in \cO $ whose reduced norm $\n (\om )$ is supported
at the prime ideals $\mathfrak p\mid \mathcal D$. The group $W_+$
may and will be regarded as the subgroup $W_+= \{ [\omega ] \in W
: \n (\omega )\in F^*_+\} $ of $W$ and both coincide whenever
$\Pic _+(F)\simeq \Pic (F)$. From the above, we obtain the exact
sequence

$$
1\rightarrow R^*_{F_+}/R_F^{*2}\rightarrow W^1\ra W_+\ra 1.
$$

Here, a totally positive unit $u\in R^*_{F_+}$ is mapped to any $\alpha _u \in \cO
^*_+$ with $\n (\alpha _u )=u$, whose existence is guaranteed
by the Hasse-Schilling-Maass Theorem in its integral version (cf.\,\cite{HaSc},
\cite{Vi}, p. 90). In this
way, we obtain that $W^1$ is isomorphic to the {\em direct} product

$$
W^1\simeq R^*_{F_+}/R_F^{*2}\times W_+\simeq (\Z /2\Z )^s
$$
for some positive integer $s\leq (n-1) + 2 r$. The bound for $s$
follows from Dirichlet's unit Theorem and the inclusion
$W_+\subseteq W$. Its precise value can actually be determined in
terms of the behaviour of signatures of elements of $R_F^*$ and
the prime ideals $\mathfrak p\mid \mathcal D$ in the class group
of $F$.

\subsection{Modular interpretation of the Atkin-Lehner group}
\label{ModularAtkin} $\\ $ In his Ph.D thesis, Jordan described
the modular interpretation of the action of the positive
Atkin-Lehner group $W_+$ on Shimura curves (\cite{Jo}). This can
be extended to the higher dimensional cases and we do so now. The
action of $B^*_+\subset \mathrm{GL }^+_2(\R )^n$ on $\mathfrak
H^n$ by Moebius transformations descends to a free action of $W^1$
on the set $\cO ^1\backslash \mathfrak H^n$ of complex points of
the Shimura variety $X_B = X_{(\cO , \cI , \varrho )}$, for any
choice of a left $\cO $-ideal $\cI $ and a positive involution
$\varrho $.

The action of an involution $\om\in W^1$ can be modularly
interpreted as follows. Let $P=[(A, \iota , \cL )]$ denote the
isomorphism class of a polarized abelian variety with
multiplication by $\cO $ viewed as a closed point of $X_{(\cO ,
\cI, \varrho )}$. Then $\om $ acts on the triplet by {\em keeping}
the same isomorphism class of the underlying abelian variety $A$
but conjugating the endomorphism structure $\iota :\cO
\hookrightarrow \End (A)$ and switching the polarization $\cL $.
Namely, $\om (P) = [(A, \iota _{\om }, \cL _{\om })]$ where

$$
\begin{matrix}
\iota _{\om }: & \cO & \hookrightarrow & \End (A) \\
    & \beta & \mapsto &
\om ^{-1} \iota (\beta ) \om
\end{matrix}
$$
and $\cL _{\om } := \frac {\om ^*(\cL )}{\n (\om )}$ is $\frac
{1}{\n (\om )}$-times the pull-back of the primitive polarization
$\cL $ by the isogeny $\om \in \cO \stackrel {\iota
}{\hookrightarrow } \End (A)$. In other words, if we regard the
first Chern class $E=c_1(\cL )$ of $\cL $ as an alternate bilinear
form on $V = \mathrm{Lie }(A(\C ))$, then $c_1(\cL _{\om }):
V\times V\longrightarrow \R $, $(u, v)\mapsto E((\om /\n (\om
))(u), \om (v))$. See Section \ref{Bpart} and \cite{Ro2} to check
that this does correspond to the first Chern class of a primitive
polarization $A$ compatible with $\iota _{\om }$. From this
interpretation and by standard moduli considerations, it follows
that $W^1\subseteq \Aut (X_B)$ acts on $X_B$ as a subgroup of
algebraic involuting automorphisms over $\Q $.

\begin{remark} If $\om \in \cO ^*_+$, $\n (\om )=u\in R^*_{F_+}$, then
$\om $ can be viewed as an automorphism $\om \in \Aut (A)\simeq \cO ^*$
of $A$. This automorphism induces an isomorphism of triplets $\om (A,
\iota , \cL ) = (A, \om ^{-1}\iota \om , \om ^*(\cL )/u) \stackrel {\om
}{\simeq }(A, \iota , \cL _u)$ where $\cL _u$ is a principal
polarization on $A$ such that, although $\cL \not \simeq \cL _u$, both
induce the same Rosati involution on $\End (A)\otimes \Q $.

In the literature (cf.\,\cite{Sh1}, \cite{Mi}), $\cL $ and $\cL _u$
are called weakly isomorphic. It is easy to see that $\cO ^*_+/\cO
^1$ acts freely and transitively on the set of isomorphism classes
of a given weak polarization class on $A$.
\end{remark}

\begin{remark} There is a natural moduli theory for
polarized abelian varieties with QM up to weak isomorphism which is also
considered in \cite{Sh1}. Both theories coincide in dimension 2 (because it
corresponds to $F=\Q $), but for higher dimensions the
latter is coarser and less suitable for our purposes.
\end{remark}

\begin{remark} We wonder in what circumstances $W^1$ is the full
group of automorphisms of $X_B$. The impression is that, {\em
generically}, it does hold that $\Aut (X_B) = W^1$, but of course
the term {\em generic} should be made precise in any case. The
split case of an order $\cO $ in $B=\mathrm{M}_2(\Q )$ was
classically studied by Ogg, who found that the modular curve
$X_0(37)$ is an interesting exception to the prediction that,
whenever the genus of a modular curve is at least 2, $W=W^1=\Aut
(X_{\cO })$. For Shimura curves $X_{\cO }$ with $\cO $ in a
rational division quaternion algebra $B/\Q $, this question was
investigated in \cite{Ro1}. For higher dimensional Shimura
varieties, it seems to hold that whenever $\cO ^*$ does not
contain torsion units (besides $\pm 1$), $W^1 = \Aut (X_{\cO , \cI
, \varrho })$ (cf.\,\cite{Ro1}, Theorem 2).
\end{remark}

\section{Main theorem and corollaries}
\label{theorems}

As remarked in the introduction, the type of the primitive
polarizations of the abelian varieties parametrized by the Shimura
variety $X_{(\cO , \cI , \varrho )}$ is determined by the datum
$(\cO , \cI , \varrho )$. The following was proved in \cite{Ro2}
and makes this observation explicit.

\begin{proposition}
\label{Ro}

The po\-lari\-zations of the abe\-lian varie\-ties with
qua\-ter\-nio\-nic mul\-tipli\-cation $(A, \iota , \cL )$
parametrized by $X_{(\cO , \cI , \varrho )}$ are {\em principal}
if and only if:

\begin{itemize}
\item
$\mathrm{Disc }(B) = (D)$ is a principal ideal of $F$ generated
by a totally positive element $D\in F^*_+$.
\item
$\n _{B/F}(\cI )$ and the codifferent
$\vartheta _{F/\Q }^{-1}=\{ x\in F, \tr _{F/\Q }(x R_F)\subseteq \Z
\}$ of $F$ over $\Q $ lie in the same ideal class in $\Pic (F)$.
\item
The positive anti-involution on $B$ is
$\varrho = \varrho _{\mu }: B\rightarrow B$, $\beta \mapsto
\mu ^{-1} \bar \beta \mu $ for some
$\mu\in \cO $ such that $\mu ^2 + D = 0$.

\end{itemize}

\end{proposition}

For the rest of the note, we will focus on moduli spaces of
principally polarized abelian varieties and therefore we place
ourselves under the assumptions on $(\cO , \cI , \varrho _{\mu })$ of
Proposition \ref{Ro}. We thus assume in particular that
$\mathrm{disc}(B)=(D)$ for some $D\in F^*_+$ and that $\mu \in \cO $ satisfies
$\mu ^2+D=0$. In this case
we say that $(\cO , \cI , \varrho _{\mu })$ is of
{\em principal type} and we will simply denote by $\mathcal A_g$ the
moduli space $\mathcal A_{g, (1, ..., 1)}$.

Our main Theorem \ref{main} below describes how the modular maps
introduced in Section 2.1 factor through the quotient of $X_B$ by
certain subgroups of $W^1$ of Atkin-Lehner involutions that we
introduce first.

Let us say that an element $\chi \in \cO \cap
\Norm _{B^*}(\cO )$ is a twist of $(\cO ,\mu )$ if
$\chi ^2+\n (\chi )=0, \mu \chi =-\chi \mu $ and therefore

$$
B = F+F\mu +F\chi +F \mu \chi = (\dfrac {-D, -\n (\chi )}{F}).
$$

Further, we will
say that a pair $(\cO , \mu )$ is {\em twisting}
if it admits some twist $\chi $ and that a quaternion algebra $B$ is {\em twisting}
if it contains a twisting pair $(\cO , \mu )$. Note that a totally
indefinite quaternion algebra is twisting if and only if
$B=(\dfrac {-D, m}{F})$ for some $m\in F^*_+$, $m\mid D$.

Now let $\om \in W^1$ be an Atkin-Lehner involution of $\cO $. We
say that it is a twisting involution with respect to $(\cO , \mu
)$ if the class of $\om $ in $W=\Norm _{B^*}(\cO )/F^* \cO ^*$ is
represented by a twist $\chi \in \cO \cap \Norm _{B^*}(\cO )$ of
$(\cO , \mu )$. For any subring $S\subset \cO $ we say that a
twist $\chi $ of $(\cO , \mu )$ is in $S$ if $\chi \in S$. An
involution $\om \in W^1$ is a twisting involution in $S$ if it can
be represented by a twist $\chi \in S$.

\begin{definition}
\begin{enumerate}
\item
The {\em twisting group} $V_0$ attached to $(\cO , \mu )$ is the subgroup of
$W^1$ generated by the twisting involutions of $(\cO , \mu )$.

\item
The {\em stable group} attached to $(\cO , \mu )$ is the subgroup

$$
W_0=V_0\cdot U_0
$$
of $W^1$ generated by
$$
U_0 = \Norm _{F(\mu )^*}(\cO )/F^*\Omega (S),
$$
where $\Omega (S)=\{ \xi \in S, \xi ^f=1, f\geq 1\} $ denotes the
finite group of roots of unity in the CM-order $S=F(\mu )\cap \cO $,
and the twisting group $V_0$.

\end{enumerate}
\end{definition}

For any subring $S\subseteq \cO $, we shall let $V_0(S)$ denote the subgroup of
$W^1$ generated by the twisting involutions of $(\cO , \mu )$ in $S$.

Let
$$
\widetilde {X}_B = \widetilde {X}_{(\cO , \cI , \varrho )} = \pi
(X_{(\cO , \cI , \varrho )})\hookrightarrow \mathcal A_g
$$
denote the image in $\mathcal A_g$ by $\pi $ of the Shimura
variety $X_{(\cO , \cI , \varrho )}$. Similarly, define
$\widetilde {X}_{B/F}\subset \mathcal H _F$ and $\widetilde
{X}_{B/(S, \varphi )}\subset \mathcal H _S$ respectively to be the
algebraic subvarieties $\pi _{F}(X_{(\cO, \cI , \varrho )})$ and
$\pi _{(S, \varphi )}(X_{(\cO, \cI , \varrho )})$ of the Hilbert
modular varieties  $\cH _{F}$ and $\cH _{S}$. Note that the
classical Humbert varieties arising from Hilbert modular varieties
can be considered as a degenerate case of the above defined
$\widetilde {X}_B$, since they would correspond to the split
algebra $B=\mathrm{M}_2(F)$.

\begin{definition}
A closed point $[(A, \iota ,\cL )]$ on $X_B$ and its image on
$\widetilde {X}_{B/(S, \varphi )}$, $\widetilde {X}_{B/F}$ or
$\widetilde {X}_B$ will be called a {\em Heegner point} if $\End
(A)\otimes \Q \simeq \M _2(L)$ for a CM-field $L$ over $F$.
\end{definition}

The set of Heegner points on these varieties is discrete and
dense, $\End (A)\otimes \Q = B$ being the generic case.

We say that a morphism between two schemes of non necessarily the
same dimension is {\em quasifinite} if it has finite fibres
(cf.\,\cite{Ha}, p.\,91).

\begin{theorem}
\label{main}

Let $(\cO , \cI , \varrho )$ be a quaternionic datum of principal
type and let $X_B= X_{(\cO , \cI , \varrho )}$ be the
corresponding Shimura variety. For any totally real Eichler pair
$(S, \varphi )$, let

$$
\begin{matrix}
     &    &        & \cH _S &  & & & \\
     &        &\stackrel {\pi _{S, \varphi }}{\nearrow } & & \searrow
     & & & \\
\pi :& X_B  & &\stackrel {\pi _F}{\longrightarrow }  & & \cH _F &
\rightarrow &\cA
_g.\\
\end{matrix}
$$
be the diagram of forgetful morphisms introduced before. Then

\begin{itemize}

\item The map  $\pi _F:X_B\rightarrow \cH _F$ is a quasifinite map that
factors over $\Q $ into the natural projection $X_B\rightarrow
X_B/W_0 $ from $X_B$ onto its quotient by the stable group
$W_0\subseteq \Aut (X_B)$ and a birational morphism $b_F:
X_B/W_0\dashrightarrow \widetilde {X}_{B/F}$ onto the image of
$X_B$ in $\cH _F$.

The maximal domain of definition of $b_F^{-1}$ is $\widetilde
{X}_{B/F}\setminus \mathcal T _F$, where $\mathcal T _F$ is a
finite set of Heegner points.

\item
The map  $\pi _{(S, \varphi )}:X_B\rightarrow \cH _S$ is a
quasifinite map that factors over $\Q $ into the projection
$X_B\rightarrow X/V_0(\varphi (S))$ of $X_B$ onto its quotient by
the finite 2-group $V_0(\varphi (S))\subseteq W_0$ and a
birational morphism $b_{(S, \varphi )}: X_B/V_0(\varphi
(S))\dashrightarrow \widetilde {X}_{B/(S, \varphi )}$ into the
image of $X_B$ in $\cH _S$ by $\pi _{(S, \varphi )}$.

As before, $b_{(S, \varphi )}^{-1}$ is defined on the whole
$\widetilde {X}_{B/(S, \varphi )}$ but at a finite set $\mathcal T
_{(S, \varphi )}$ of Heegner points.

\end{itemize}

\end{theorem}

We call the sets $\mathcal T _F$ and $\mathcal T _{(S, \varphi )}$
the {\em singular Heegner loci} of $\widetilde {X}_{B/F}$ and
$\widetilde {X}_{B/(S, \varphi )}$ respectively, as they are
indeed sets of singular points (of quotient type) on these
varieties. We present the proof of the theorem in Sections
\ref{Bpart} and \ref{proofs}; we now derive some corollaries of
it.

Being subgroups of $W^1$, both the twisting subgroups $V_0(\varphi
(S))$ and the stable subgroup $W_0$ are isomorphic to the direct
product of a certain number of copies of $C_2$. We refer the
reader to \cite{Ro3} for a detailed description of them. As a
consequence of the results proved in \cite{Ro3} and our
main Theorem \ref{main}, we obtain the
following.

\begin{corollary}

Let $(\cO , \mu)$ be a non twisting polarized order. Then, the
maps $\pi _{(S, \varphi )}:X_B\rightarrow \cH _S$ are birational
equivalences for any totally real quadratic order embedded in $\cO
$ and

$$
\deg (\pi _F:X_B\rightarrow \cH _F) = 2^{\omega _{odd}},
$$
where $\omega _{odd} = \sharp \{ \xi \in R_{\mu }, \xi ^{f}=1,
f\mbox{ odd }\} $.

$\\ $ Let $(\cO , \mu)$ be a twisting polarized order and let
$\chi \in \cO $ be any twist of $(\cO , \mu )$. Then, $\pi
_{R_F[\chi ]}: X_B\rightarrow \cH _{R_F[\chi ]}$ is a $2:1$ map,
whereas for any totally real Eichler pair $(S, \varphi )$ such
that $\varphi (S)$ does not contain any twist, $\pi _{(S, \varphi
)}$ is birational.

Moreover,
$$
\deg (\pi _F:X_B\rightarrow \cH _F) = 2^{2 \omega _{odd}}.
$$

\end{corollary}

In particular, note that the forgetful map $\pi _{\Q }:
X_B\longrightarrow \cA _2$ from a Shimura curve into Igusa's
moduli space of principally polarized abelian surfaces is either
of degree $2$ or of degree $4$ and that a necessary condition for
the latter is that $B\simeq (\frac{-\disc (B), m}{\Q })$ for some
$m>0$, $m\mid \disc (B)$. It turns out, for instance, that for any
choice of a datum $(\cO , \cI , \varrho )$ of principal type of
discriminant $D=6$ or $10$, the corresponding forgetful map of the
Shimura curve $X_{(\cO , \cI , \varrho )}$ into $\cA _2$ has
degree $4$.

\section{The N\'{e}ron-Severi group of an abelian variety with
quaternionic multiplication}
\label{Bpart}

In this section we introduce some ingredients that may be useful
to make the forgetful maps and the singular Heegner loci $\mathcal
T _F$ and $\mathcal T _{(S, \varphi )}$ of Theorem \ref{main}
explicit in particular cases, although we will not pursue this
purpose in this note. The content of this section will also be
used in the proof of Theorem \ref{main} in Section \ref{proofs}.

For an abelian variety $A$, the N\'{e}ron-Severi group $\mathrm{NS}(A)$ is
the group $\mathrm{Pic}(A)/\mathrm{Pic}^0(A)$ of invertible sheaves on
$A$ up to algebraic equivalence. By N\'{e}ron's Basis Theorem,
$\mathrm{NS}(A)$ is a $\Z $-module of finite rank and we denote $\NS
^0(A) = \NS (A)\otimes _{\Z }\Q $.

As always, let $\cO $ denote a maximal order in a totally indefinite
quaternion algebra $B$ over a totally real number field $F$, $[F:\Q
]=n$. Assume that $A$ is an abelian variety over $\C $ of dimension $2 n$ and
that $\End (A)$ contains an order isomorphic to $\cO $. Then $A$ is
either simple or isogenous to the square $A_0^2$ of a simple abelian
variety $A_0$ of dimension $n$. In
the former case we actually have that $\End (A)\simeq \cO $ while in
the latter $\End (A)$ is an order in a CM-field $L$ over $F$ and $\End
(A)$ can be identified with an order in $\M _2(L)$.

More precisely, and in a unified way, fix a ring homomorphism
$\iota :\cO \hookrightarrow \End (A)$ and let $R = \End _{\cO }(A)
= \{ \gamma \in \End (A): \gamma \beta = \beta \gamma \mbox{ for
all } \beta \in \iota (\cO )\}$ be the commutator of $\cO $ in
$\End (A)$. Then either $R = R_F$ or $R$ is an order in $L$ but in
any case $\End (A) = \cO \otimes _{R_F} R$.

We now wish to select a suitable piece of the N\'{e}ron-Severi group
of $A$ and describe it also in a unified way, regardless of
whether $A$ is simple or isogenous to the square of a CM abelian
variety. Recall that a line bundle $\cL \in \NS ^0(A)$ induces an
involution $\circ _{\cL } $ on $\End ^0(A)$, called the Rosati
involution: $\circ _{\cL }: \End ^0(A) \rightarrow \End ^0(A) $,
$\beta \mapsto \beta ^{\circ } = \varphi ^{-1}_{\cL } \hat {\beta
}\varphi _{\cL }$. Here, if $\hat A$ is the dual abelian variety
of $A$, we let $\hat \beta $ denote the dual endomorphism $\hat
\beta :\hat A\rightarrow \hat A$ of $\beta $ on $\hat A$ and we
let $\varphi _{\cL }: A \rightarrow \hat A $, $x \mapsto \tau
_x^*(\cL )\otimes {\cL }^{-1}$ be the morphism that $\cL $ induces
from $A$ to $\hat A$. By $\tau _x$ we mean the translation-by-$x$
map on $A$.

\begin{definition}

Let $(A, \iota )$ be an abelian variety of dimension $2 n$ endowed with
a multiplication by $\cO $. The {\em $B$-part of the (rational)
N\'{e}ron-Severi group of $A$} is

$$
\NS _B^0(A) = \{ \cL  \in \NS ^0(A), \circ _{\cL } (B) = B\} .
$$
The
{\em $B$-part of the N\'{e}ron-Severi group of $A$} is the group $\NS _B(A) = \NS
^0_B(A)\cap \NS (A)$ of line bundles $\cL \in \NS (A)$ on
$A$ whose Rosati involution $\circ _{\cL }$ leaves $B$ invariant.

\end{definition}

If $A$ is simple (and thus $\iota :\End (A)\simeq \cO $ is an
isomorphism), we do not obtain anything new: $\NS _B(A) = \NS (A)$
and $\NS ^0_B(A) = \NS ^0(A) $. However, when $A$ is not simple
(and thus $\End ^0(A)\simeq \M _2(L)$ for some CM-field $L$),
these are proper subgroups of the usual N\'{e}ron-Severi groups. The
following extends a result proved by the author in \cite{Ro2} and
allows us to see line bundles on $A$ as pure quaternions in $B$ in
a canonical way. See Section 1 for notations.

\begin{proposition}
\label{NS}

Let $(A, \iota )$ be a complex abelian variety of dimension $2 n$
together with a multiplication by $\cO$. Let $\cI $ be the isomorphism
class of the left $\cO $-module $H_1(A, \Z )$, represented by a left
$\cO $-ideal $\cI \subset B$. Then the
first Chern class induces a natural isomorphism of groups

$$
\NS _B (A)\stackrel {\simeq }{\rightarrow }\quad \cN (\cI )^{\sharp
}_0\hookrightarrow \quad B_0.
$$

\end{proposition}

Explicitly, the isomorphism is constructed as follows: write
$A(\C) = V/\Lambda $ for a vector space $V$ and a lattice $\Lambda
= \cI v_0$, $v_0\in V$, and let $\cL \in \NS _B(A)$. Let also
$E=c_1(\cL ):V\times V\rightarrow \R $ be its first Chern class
regarded as an $\R $-alternate Riemann form on $V$. Since
$E(\Lambda \times \Lambda )\subseteq \Z$, we may consider the
linear map $B\rightarrow \Q $ that sends an element $\beta \in B$
to $E(\iota (\beta ) v_0, v_0)$. By the non-degeneracy of the
reduced trace $\tr _{B/\Q } = \tr _{F/\Q }\cdot \tr _{B/F}$, we
can find an element $\mu \in B$ such that $E(\iota (\beta ) v_0,
v_0) = \tr _{B/\Q }(\mu \beta )$ for any $\beta \in B$. It follows
from the alternateness of $E$ that $\mu $ is pure: $\mu \in B_0$
(cf.\,\cite{Ro2} for details).

Since $\cL \in \NS _B(A)$, the Rosati involution $\circ _{\cL }$ on $\End
^0(A)$ descends to a well-defined involution $\circ :B\rightarrow B$
such that, for any $u$, $v\in V$, $E(u, \iota (\beta ) v) = E(\iota
(\beta ^{\circ }) u, v)$. By means of the Skolem-Noether Theorem it can
be seen that $\beta ^{\circ } = \mu ^{-1} \bar \beta \mu $ (cf.\,\cite{Ro2}
again for details). This already determines the Riemann form $E=c_1(\cL )$:

$$
\begin{matrix}
E: & V\times V & \longrightarrow & \R \\
   & (u, v) & \mapsto & \tr _{B\otimes _{\Q }\R /\R }(\mu \gamma \bar \beta )
\end{matrix}
$$
where $\beta$ and $\gamma \in B\otimes _{\Q }\R \simeq \M _2(\R
)^n$ are uniquely determined elements such that $\beta v_0=u$ and
$\gamma v_0 = v$. Observe that, since $\Lambda = \cI v_0$ and
$E(\Lambda \times \Lambda )\subseteq \Z $, $\mu \in \cN (\cI
)^{\sharp}_0$.

In what follows, the above Riemann form will be denoted as $E_{\mu
}$. By an abuse of notation, we will also sometimes use the
notation $c_1(\cL )=\mu $. It was shown in \cite{Ro2}, Theorem 2.2, that the
assignation $\cL \mapsto c_1 (\cL )$ is an isomorphism of
groups between $\NS _B(A)$ and $\cN (\cI )^{\sharp }_0$ in the
case that $A$ is simple. As the reader may check, the proof is
exactly the same as in the CM split case.

Since the rank of $\NS _B(A)$ over $\Z $ is $3 [F:\Q ]= 3 n$, {\em
any }$\Z$-module of the same rank is actually isomorphic to it.
The interest of the above description of the $B$-part of the
N\'{e}ron-Severi group $\NS _B(A)$ is that the isomorphism of
Proposition \ref{NS} behaves nicely with respect to the degree of
the line bundles $\cL \in \NS _B(A)$ and also with respect to the
pull-back $\alpha ^*(\cL )$ by endomorphisms $\alpha \in \End (A)$
of $A$. Indeed, let us define the degree $\mathrm{deg }(\cL )$ of
a non-degenerate line bundle $\cL \in \NS (A)$ to be the degree of
the finite map $\varphi _{\cL }: A \rightarrow \hat A$ from $A$ to
the dual abelian variety $\hat A$ induced by $\cL $.

\begin{proposition}
\label{deg}

If the first Chern class of $\cL $ is represented by the Riemann
form $E_{\mu }$ with $\mu \in B_0$ and $\mu ^2+\delta =0$, $\delta
\in F^*$, then $\mathrm{deg }(\varphi _{\cL })={\mathrm{N}_{F/\Q
}( \vartheta _{F/\Q }^2\cdot \mathrm{n}_{B/F}(\cI )^2\cdot
\mathcal D\cdot \delta )}^2$.

\end{proposition}

The proof of Proposition \ref{deg} was given in \cite{Ro2}, Proposition
3.1, for non-CM abelian varieties but the same arguments are valid for
this case as well. In order to describe now how the first Chern
class behaves under pull-backs of line bundles by endomorphisms,
we denote by $\lambda \mapsto \bar {\lambda }\in L$ the complex
conjugation on $L$. There should be no confusion with the
quaternionic conjugation $\beta \mapsto \bar {\beta }$ on $B$.
Both conjugation maps combine to endow $B\otimes _F L \simeq \M
_2(L)$ with a conjugation map $\lambda \beta \mapsto \bar {\lambda
}\bar {\beta }$ (and extended linearly). In matrices, $\overline
{\begin{pmatrix}
  a & b \\
  c & d
\end{pmatrix}} = \begin{pmatrix}
  \bar d & -\bar b \\
  -\bar c & \bar a
\end{pmatrix}\in \M _2(L)$.

Having fixed a homomorphism of algebras $\iota : B \hookrightarrow \End
^0(A)$, identify $\End ^0(A)$ with either $B$ or $\M _2(L)$ with $B$ as
an embedded sub-algebra.

\begin{proposition}
\label{pullbacks}

Let $(A, \iota )$ be an abelian variety with quaternionic
multiplication by a maximal order in $B$ and let $\cL \in \NS
_B^0(A)$ be a line bundle on $A$ compatible  with $B$. Denote by
$c_1(\cL ) = \mu \in B_0$ its first Chern class regarded as a pure
quaternion in $B$.

Then the pull-back polarization $\alpha ^*(\cL )$ of $\cL $ by a
rational endomorphism $\alpha \in \End ^0(A)$ is compatible with
$B$, that is, $\alpha ^*(\cL )\in \NS _B^0(A)$, if and only if
$\bar {\alpha }\mu \alpha \in B_0$. In this case, $c_1(\alpha
^*(\cL )) = \bar {\alpha }\mu \alpha $.

\end{proposition}

Obviously, when $A$ is simple the above condition is empty and all
polarizations are compatible with $B$. Moreover, in the simple
case the relation $c_1(\alpha ^*(\cL )) = \bar {\alpha }c_1(\cL )\alpha
$ was already proved in \cite{Ro2}, Theorem 2.2. However, the proof of this
fact in the CM case is not trivial and deserves a more careful
inspection.

We hence assume that $\End ^0(A)=\M _2(L)$ for a CM-field $L$ over
$F$ and we first claim that $\alpha ^*(\cL )$ induces the Rosati
involution $\beta \mapsto (\bar {\alpha }\mu \alpha )^{-1} \bar
{\beta } (\bar {\alpha }\mu \alpha )$ on $\M _2(L)$. Denote
$E_{\cL }$ and $E_{\alpha ^*(\cL )}$ the Riemann forms on $V =
\mathrm {Lie} (A)$ induced by $\cL $ and $\alpha ^*(\cL )$
respectively. The Rosati involution $\circ _{\alpha ^*(\cL )}$
restricts to complex conjugation on $L$ and is characterized by
the fact that $E_{\alpha ^*(\cL )}(u, \beta v) = E_{\alpha ^*(\cL
)}(\beta ^{\circ _{\alpha ^*(\cL )}} u, v)$ for any $u$, $v\in V$
and $\beta \in \End ^0(A)$.

The above given anti-invo\-lu\-tion is in\-deed com\-plex
con\-ju\-ga\-tion on $L$ and it holds that $E_{\alpha ^*(\cL )}(u,
\beta v) = E_{\cL }(\alpha u, \alpha  \beta v)$ $= E_{\cL }(\mu
^{-1} \bar {\beta } \bar {\alpha } \mu  \alpha u, v) $ $=
E_{\alpha ^*(\cL )}((\bar {\alpha } \mu  \alpha )^{-1} \bar {\beta
} (\bar {\alpha } \mu  \alpha )u, v)$ for any $u$, $v\in V$ and
$\beta \in \M _2(L)$.

Now that we know this, $\alpha ^*(\cL )\in \NS _B^0(A)$ if and
only if $\bar {\alpha } \mu \alpha \in \M _2 (L)$ normalizes $B$.
But the Skolem-Noether Theorem implies that an element in $\M
_2(L)$ normalizes $B$ if and only if it belongs to the
multiplicative group $L^*\cdot B^*$.

In order to see that actually $\bar {\alpha } \mu \alpha \in B$
(and thus also in $B_0$), write $L = F( \eta )$, $\eta ^2 =
-\Delta \in F^*_+$, and express $\alpha = \alpha _0 + \alpha _1
\eta \in B\otimes _F L \simeq \M _2(L)$ with $\alpha _i\in B$.
Then the reader may check that $\bar {\alpha } \mu \alpha  = (\bar
{\alpha _0} \mu \alpha _0 + \Delta \bar {\alpha _1} \mu \alpha _1)
+ \tr _{B/F}(\bar {\alpha _0} \mu \alpha _1)\eta $. Since $\bar
{\alpha } \mu \alpha $ must belong to $L^*\cdot B^*$, then either
$\bar {\alpha } \mu \alpha = (\bar {\alpha _0} \mu \alpha _0 +
\Delta \bar {\alpha _1} \mu \alpha _1)$ or $\bar {\alpha } \mu
\alpha = \tr _{B/F}(\bar {\alpha _0} \mu \alpha _1)\eta $.

We claim that the latter cannot occur and therefore $\bar {\alpha
} \mu \alpha \in B_0$, as we wish. For if $\bar {\alpha _0} \mu
\alpha _0 = -\Delta \bar {\alpha _1} \mu \alpha _1$, this would
imply that

\begin{itemize}
\item
$\n _{B/F}(\alpha _0) = \Delta \n _{B/F}(\alpha _1)$
and
\item
$\mu \alpha _o \alpha _1^{-1} = -\Delta \bar {\alpha _0}^{-1}\alpha _1
\mu = -\alpha _0 \alpha _1^{-1} \mu $.
\end{itemize}

It would then follow that $\tr _{B/F}(\alpha _0 \alpha _1^{-1}) =
0$ and hence $B = (\frac {-\Delta , -\Delta }{F})$, which is a
con\-tra\-dic\-tion be\-cause $B$ is to\-tally indefinite whi\-le
$(\frac {-\sigma (\Delta ), -\sigma (\Delta )}{\R })$ is the skew
field of Hamilton quaternions, for any real embedding $\sigma
:F\hookrightarrow \R $.

This proves the first part of the proposition. Assume now that
$\alpha ^*(\cL )\in \NS _B^0(A)$ and thus $\bar {\alpha } \mu
\alpha  = \bar {\alpha _0} \mu \alpha _0 + \Delta \bar {\alpha _1}
\mu \alpha _1\in B_0$. In order to show that indeed $c_1(\alpha
^*(\cL )) = \bar {\alpha } \mu \alpha $, recall that for a
polarization $\cL \in \NS _B^0(A)$, the Riemann form associated to
it is $E_{\mu }: V\times V\longrightarrow \R $, $(u, v) \mapsto
\tr _{B/\Q }(\mu \gamma \bar {\beta })$ where $\mu = c_1(\cL )$
and $\beta $, $\gamma \in B\otimes _{\Q }\R $ are the only
elements such that $\beta v_0 = u$ and $\gamma v_0 = v$
respectively.

We must show that $E_{\bar {\alpha } \mu \alpha } (u, v) = E_{\mu
}(\alpha u, \alpha v)$ for any $u$, $v\in V$. This follows from
direct computation provided one takes into account that, if we let
$\zeta \in B\otimes _{\Q }\R $ denote the only element such that
$\zeta v_0 = \eta v_0$, then for any $\beta \in B\otimes _{\Q }\R
$ we have that $\alpha \beta v_0 = (\alpha _0 \beta + \alpha _1
\beta \zeta ) v_0$.

\section{Proof of Theorem \ref{main}}
\label{proofs}

Let $(\cO , \cI , \varrho )$ be a quaternionic datum attached to a
totally indefinite quaternion algebra $B$ over a totally real
number field $F$. We assume that it is of principal type
(cf.\,Proposition \ref{Ro} and below). This means in particular that
$\mathrm {disc }(B) = (D)$ can be generated by an element $D\in F^*_+$
such that $\varrho = \varrho _{\mu }$ for some $\mu \in \cO $, $\mu
^2+D=0$.

Let $X_B = X_{(\cO , \cI , \varrho )}$ be the Shimura variety
attached to it. By Section \ref{Atkin}, the elements in the
totally positive Atkin-Lehner group $W^1$ act as a group of
rational automorphisms on $X_B$. Let us first show that $\pi
_F:X_B\rightarrow \cH _F$ is a quasifinite map that factors into
the natural projection

$$
X_B\rightarrow X_B/W_0
$$
from $X_B$ onto its quotient by the stable group $W_0 = U_0\cdot
V_0 \subseteq \Aut _{\Q }(X_B)$ attached to $(\cO , \varrho )$.

To this end, we  claim that there is an action of the stable group
$W_0$ on the geometric fibres of the morphism $\pi _F:
X_B\rightarrow \cH _F$ that is free and transitive on a certain
open Zariski-dense subset $\mathcal U _F$ of $X_B$.

In order to prove the claim, let $(A, j_F, \cL )$ be a principally
polarized abelian variety together with a homomorphism of $\Z
$-algebras $j_F:R_F\hookrightarrow \End (A)$. Without making any
further mention of it, we will regard $\End (A)$ as an
$R_F$-algebra through the given immersion $j_F$.

The isomorphism class $[(A, j_F, \cL )]$ may be interpreted as a
closed point in $\cH _F$. If it is non-empty, the elements in
$X_B$ of the fibre of $\pi _F$ at this point can then be
interpreted as those principally polarized abelian varieties with
quaternionic multiplication whose isomorphism class is represented
by a triplet $(A, \iota , \cL )$ where $\iota : \cO
\hookrightarrow \End (A)$ is a homomorphism of $R_F$-algebras such
that the Rosati involution that $\cL $ induces on $\cO $, via
$\iota $, coincides with $\varrho _{\mu }$.

Choose a triplet $(A, \iota , \cL )$ as above. Proposition
\ref{NS} identifies $\NS (A)$ with a lattice in $B_0$ by means of a map that
we have called $c_1$. We warn that this identification depends on $\iota $ and
it allows us to regard polarizations on the abelian variety
$A$ as pure quaternions of $B$. In particular, since $\circ
|_{\iota (\cO )}=\varrho _{\mu }\cdot \iota $, we have that
$c_1(\cL )=\mu $ up to multiplication by elements in $F^*$.

Recall now that $W_0 = U_0\cdot V_0$ for certain subgroups $U_0$,
$V_0\subseteq W^1$ that were defined in Section \ref{Atkin}. Let
$\omega \in U_0$. Let $L = F(\mu )\simeq F(\sqrt {-D})$ be the
CM-field generated by $\mu \in B$ over $F$ and let $S\supseteq
R_F[\mu ]$ be the order in $L$ at which $\mu $ is optimal. Then
$\omega $ is represented by an element (that we still denote)
$\omega \in S$ and we wish to show that the closed points $[(A,
\iota , \cL )]$ and $\om [(A, \iota , \cL )]$ in $X_B$ lie on the
same fibre of $\pi _F$.

We re\-call (cf.\,Section \ref{ModularAtkin}) that the iso\-mor\-phism class
of $\om [(A, \iota , \cL )]$ is
represented by a triplet $(A, \iota _{\om }, \cL _{\om })$ where
$\iota _{\om } = \om ^{-1}\iota \om :\cO \hookrightarrow \End (A)$ and
$\cL _{\om } =
\frac {\om ^*(\cL )}{\n (\om )}$ is a principal polarization on $A$.
By Proposition \ref{pullbacks}, it is such that $c_1(\cL _{\om }) =
\frac {\bar {\om }}{\n (\om )} c_1(\cL ) \om  =
\om ^{-1} c_1(\cL ) \om = c_1(\cL ) $.
Here, the last equality holds because $c_1(\cL )=\mu \in B^*/F^^*$ and $\mu $ and
$\om $ belong to the same quadratic order $S$ embedded in $\cO $ and hence commute.
Therefore $c_1(\cL ) = c_1(\cL _{\om })$ and we conclude
that $(A, \cL )$ and $(A, \cL _{\om
})$ are isomorphic polarized varieties. Moreover, $\iota $ and $\iota _{\om
}$ coincide when restricted to the centre $R_F$ of $\cO $ and we obtain that
$(A, \iota |_{R_F}, \cL )\simeq (A, \iota _{\om }|_{R_F}, \cL _{\om })$, which is
what we needed to show.

Now let $\omega \in V_0$ be an element represented by $\omega \in
\cO _+$ with $\n (\om ) = m\in F^*_+$ and such that $B=F+F \mu +F
\chi +F \mu \chi $ for some $\chi \in \cO $, $\chi ^2=m$ , $\mu
\chi = -\chi \mu $. Let $[(A, \iota , \cL )]\in X_B$ be a closed
point over $[(A, j_F, \cL )]\in \cH _F$. We have that $\om [(A,
\iota , \cL )] = [(A, \iota _{\om }, \cL _{\om })]$ and we must
show that $(A, \iota |_{R_F}, \cL )\simeq (A, \iota _{\om
}|_{R_F}, \cL _{\om })$. Again, since $\iota |_{R_F}=\iota _{\om
}|_{R_F}$, we only need to see that $\cL $ and $\cL _{\om }$ are
isomorphic polarizations on $A$. To do so, we first note that
$\alpha = \om ^{-1}\cdot \chi $ is a unit in $\cO $ of reduced
norm $\n (\alpha )=-1$. This is due to the fact that the principal
ideal $m R_F$ is supported at the prime ideals $\mathfrak p\mid
D$, as can be checked locally. This unit, or rather $\iota (\alpha
)$, is an automorphism of $A$ and we now show that $\alpha ^*\cL
_{\omega }=\cL \in \NS (A)$. Indeed, by Proposition
\ref{pullbacks} and since $c_1(\cL )$ and $\chi $ anti-commute,
$c_1(\alpha ^*(\cL _{\om })) = \bar {\alpha } \om ^{-1} c_1(\cL )
\om \alpha =\frac {1}{m}\bar {\chi } c_1(\cL ) \chi  = -\chi ^{-1}
c_1(\cL ) \chi  = \chi ^{-1} \chi  c_1(\cL ) = c_1(\cL )$. Since
$W_0=U_0\cdot V_0$, this proves part of our claim.

In order to show that the action of $W_0$ on the fibres of $\pi
_F$ is free and transitive in an open and dense subset $\mathcal U
_F$ of $\widetilde {X}_{B/F}\subset \cH _F$, assume that $[(A,
j_F, \cL )]$ is a point in $\widetilde {X}_{B/F}$ that is
represented by an abelian variety $A$ such that $\iota :\cO
\stackrel {\simeq }{\rightarrow }\End (A)$ is an isomorphism. In
other words, we consider the fibre of $\pi _F$ at a non-Heegner
point. In this case, if $[(A, \iota , \cL )]$ lies over $[(A, j_F,
\cL )]$, then, by the Skolem-Noether Theorem, any element on the
fibre of $\pi _F$ at this point must be represented by a triplet
$(A, \om ^{-1}\iota \om , \cL _{\om })$ for some $\om \in
\mathrm{Norm} _{B^*}(\cO )$. Moreover, for $\cL _{\om }$ to be a
polarization on $A$, $\om $ must have totally positive reduced
norm (cf.\,\cite{Ro2}, Theorem 5.1). Since we clearly have that
$(A, \iota , \cL )\simeq (A, \iota _{\om }, \cL _{\om })$ for any
$\om \in F^* \cO ^1$, we deduce that there exists a subgroup of
$W^1$ that acts freely and transitively on the fibre of $\pi _F$
at any non-Heegner point.

Therefore, this subgroup must contain $W_0$. Let us now see that
it cannot be larger than $W_0$. Assume that $\om \in W^1$ is such
that $\cL \simeq \cL _{\om }$. Then there exists $\alpha \in \cO
^*\stackrel {\iota }{\simeq }\Aut (A)$ such that $\bar {\alpha
}\om ^{-1}\mu \om \alpha =\mu $. Taking reduced norms, this
already implies that $\n (\alpha )^2=1$.

If $\n (\alpha )=1$, then $\bar {\alpha }=\alpha ^{-1}$ and the
above yields $\mu \om \alpha =\om \alpha \mu $. This means that
$\om \alpha \in F(\mu )\cap \cO = S$ and therefore $\om \in U_0$.
If $\n (\alpha )=-1$, then $\mu \om \alpha = -\om \alpha \mu $.
Write $\chi = \om \alpha $. From the discussion above, we have
that $\tr (\mu \chi )= \mu \chi + \bar {\chi }\bar {\mu }= \mu
\chi -\bar {\chi }\mu = -\tr (\chi )\mu \in F$. Thus, since $\tr
(\chi )\in F$ as well, we deduce that actually $\tr (\chi )=0$.
Since $\n (\chi )=\n (\om ) \n(\alpha ) = -\n (\om )$, we obtain
that $\chi ^2=\n(\om )$. This says that $\omega \in V_0$.

Now that we know this, we obtain as a consequence that $\pi _F:
X_B \rightarrow \cH _F$ is a quasifinite map that factors through
the natural projection of $X_B$ onto the quotient $X_B/W_0$ and a
morphism $b_F:X_B/W_0 \rightarrow \cH _F$ that is one-to-one
outside the Heegner locus of $X_B/W_0$. It is a well-known fact
(see \cite{Sh2}) that Heegner points on $X_B$ (and thus also on
the quotient $X_B/W_0$) are isolated. Since $b_F$ is an algebraic
morphism, it must be a birational equivalence between $X_B/W_0$
and its image in $\cH _F$ whose inverse must be defined everywhere
but at a finite set $\mathcal T _F$ of Heegner points.

Moreover, since $W_0\subseteq \Aut _{\Q }(\mathcal S _B)$, the
projection $\mathcal S _B\rightarrow \mathcal S _B/W_0$ is defined
over $\Q $. Since $\pi _F$ is also a morphism over $\Q $ which is
the composition of the above projection and the birational
equivalence $b_F$, it follows that $b_F$ is also defined over $\Q
$. This finishes the proof of the first part of Theorem
\ref{main}.

Now let $(S, \varphi )$ be an Eichler pair for $\cO $ and identify
$S$ with its image $\varphi (S)$ in $\cO $. As we saw, it induces
a natural morphism $\pi _{S, \varphi }: X_B\rightarrow \cH _S$
from $X_B$ into the Hilbert modular variety $\cH _S$ in such a way
that we have $\pi _F: X_B\stackrel {\pi _{(S, \varphi )}}
{\rightarrow }\cH _S\rightarrow \cH _F$. As the situation is very
similar to the one studied above, we will limit ourselves to
showing that the subgroup $V_0(S)$ of the stable group $W_0$ acts
freely and transitively on the fibre of $\pi _{(S, \varphi )}$ at
any non-Heegner point of $\pi _{(S, \varphi )}(X_B) = \widetilde
{X}_{B/(L, \varphi )}$ in $\cH _L$. This already implies that $\pi
_{(S, \varphi )}$ is quasifinite.

Let $[(A, \iota , \cL )]$ be a closed point on $X_B$ represented
by a principally polarized abelian variety with quaternionic
multiplication. From the above discussion, it is clear that any
other point at the same geometric fibre of $\pi _{(S, \varphi )}$
as $[(A, \iota , \cL )]$ is represented by $(A, \iota _{\om }, \cL
_{\om })$ for some $\om \in W_0$.

If $(A, \iota \cdot \varphi , \cL )\simeq (A, \iota _{\om }\cdot \varphi , \cL _{\om
})$ is an isomorphism given by $\alpha \in \cO ^*$, then, as we already
saw, $\n (\alpha )=\pm 1$. Write $\chi =\om \alpha $.

If $\n (\alpha )=1$, $\chi \in S\subset F(\mu )$. At the same
time, we must have that, for any $\beta \in S\subset \cO $, $\chi
\beta = \beta \chi $. This means that $\chi $ commutes
element-wise with $L$ and $F(\mu )$. Since the two quadratic
extensions are distinct because the first is totally real while
the second purely imaginary, $\chi \in R_F^*$ and thus $\om $ is
the identity element of $W_0$.

If $\n (\alpha )=-1$, then $\alpha ^*(\cL_{\om })=\cL $ implies
that $B=F+F\mu +F\chi +F\mu \chi $ while $\alpha ^{-1}\iota _{\om
}|_{S} \alpha = \iota |_{S}$ says that $\chi \in S$. Thus $\om \in
V_0(S)$. The converse also holds and the theorem follows as
before. $\Box $


\begin{thebibliography}{99}

\bibitem[BaBo]{BaBo}
W.\ L.\ Baily, A.\ Borel, Compactification of arithmetic quotients of
bounded symmetric domains {\em Ann.\ Math.\ } {\bf 84} (1966), 443-507.


\bibitem[ChFr1]{ChFr1}
T.\ Chinburg, E.\ Friedman, The smallest arithmetic hyperbolic
3-orbifold, {\em Invent.\ Math.\ } {\bf 86} (1986), 507-527.

\bibitem[ChFr2]{ChFr2}
T.\ Chinburg, E.\ Friedman, Hilbert symbols, class groups and
quaternion algebras, {\em J.\ Th\'{e}or.\ Nombres Bordeaux} {\bf 12}
(2000), 367-377.

\bibitem[Cl]{Cl}
P.L.\ Clark, {\em Local and global points on moduli spaces of
abelian surfaces with potential quaternionic multiplication},
Harvard PhD.\ Thesis, in progress.


\bibitem[Ei1]{Ei1}
M.\ Eichler, Bestimmung der Idealklassenzahl in gewissen normalen
einfachen Algebren, {\em J.\ reine angew.\ Math.\ } {\bf 176} (1937),
192-202.

\bibitem[Ei2]{Ei2}
M.\ Eichler, \"{U}ber die Idealklassenzahl hypercomplexer Systeme, {\em
Math.\ Z.\ } {\bf 43} (1938), 481-494.


\bibitem[vdGe]{vdGe}
G.\ van der Geer, {\em Hilbert modular surfaces}, Ergebn.\ Math.\
Grenz. {\bf 16}, Springer 1987.

\bibitem[Ha]{Ha}
R.\ Hartshone, {\em Algebraic Geometry,} Grad.\ Texts Math.\ {\bf
52}, Springer, 1977.

\bibitem[HaHaMo]{HaHaMo}
Y.\ Hasegawa, K.\ Hashimoto, F.\ Momose, Modular conjecture for $\Q
$-curves and QM-surfaces, {\em International J.\ Math.} {\bf
10} (1999), 1011-1036.

\bibitem[HaSc]{HaSc}
H.\ Hasse, O.\ Schilling, Die Normen aus einer normalen
Divisionsalgebra, {\em J.\ reine angew.\ Math.\ } {\bf 174} (1936),
248-252.

\bibitem[Jo]{Jo}
B.W.\ Jordan, {\em On the Diophantine arithmetic of Shimura curves},
Harvard PhD.\ Thesis, 1981.


\bibitem[LaBi]{LaBi}
H.\ Lange, Ch.\ Birkenhake, {\em Complex Abelian Varieties},
Grundl. math.\ Wiss.\ {\bf 302}, Springer, 1992.

\bibitem[Mi]{Mi}
J.\ S.\ Milne, Points on Shimura varieties mod p,
{\em Proc.\ Symp.\ Pure Math.} {\bf 33} (1979), 165-184.

\bibitem[Ri]{Ri}
K.\ A.\ Ribet, , On modular representations of $\Gal (\bar \Q /\Q )$
arising from modular forms, {\em Invent.\ Math.\ }
{\bf 100} (1990), 431-476.

\bibitem[Ro1]{Ro1}
V.\ Rotger, On the group of automorphisms of Shimura curves and
applications, {\em Compositio Math.\ }, {\bf 132} (2002), 229-241.

\bibitem[Ro2]{Ro2}
V.\ Rotger, Quaternions, polarizations and class numbers, {\em J.\
reine angew.\ Math.\ } {\bf 561}, 2003.

\bibitem[Ro3]{Ro3}
V.\ Rotger, The field of moduli of quaternionic multiplication on
abelian varieties, {\em International J.\ Math.\
M.\ Sc.\ } {\bf 52} (2004), 2795-2808.

\bibitem[Ro4]{Ro4}
V.\ Rotger, Shimura curves embedded in Igusa's threefold, preprint
2002, {\em Modular curves and Abelian varieties}, Progress in
Mathematics, Birkh\"{a}user, {\bf 224} (2003), 263-273.



\bibitem[Sh1]{Sh1}
G.\ Shimura, On analytic families of polarized abelian varieties and
automorphic functions, {\em Ann.\ Math.\ } {\bf 78} (1963), 149-192.

\bibitem[Sh2]{Sh2}
G.\ Shimura, Construction of class fields and zeta functions of
algebraic curves, {\em Ann.\ Math.\ } {\bf 85} (1967), 58-159.


\bibitem[Vi]{Vi}
M.F.\ Vign\'{e}ras, {\em Arithm\'{e}tique des alg\`{e}bres de quaternions},
Lect.\ Notes Math.\ {\bf 800}, Springer, 1980.


\end{thebibliography}
\end{document}